\documentclass[a4paper]{article}
\newcommand{\proof}{\noindent {\bf Proof: }}

\newtheorem{theorem}{Theorem}

\newtheorem{lemma}{Lemma}
\newtheorem{defi}{Definition}

\def\qed{\hfill $\Box$}

\usepackage[]{amssymb}
\usepackage{graphicx}

\begin{document}
\title{Malfatti's problem on the hyperbolic plane\footnote{Dedicated to the memory of my colleague, patron and friend Istv\'an Reiman.}}
\author{\'Akos G.Horv\'ath }
\date{March 11, 2013}

\maketitle

\begin{abstract}
More than two centuries ago Malfatti (see \cite{malfatti}) raised and solved the following problem (the so-called Malfatti's construction problem):Construct three circles into a triangle so that each of them touches the two others from outside moreover touches two sides of the triangle too.
It is an interesting fact that nobody investigated this problem on the hyperbolic plane, while the case of the sphere was solved simultaneously with the Euclidean case. In order to compensate this shortage we solve the following exercise: {\em Determine three cycles of the hyperbolic plane so that each of them touches the two others moreover touches two of three given cycles of the hyperbolic plane.}
\end{abstract}

{\bf MSC(2000):} 51M10, 51M15

{\bf Keywords:} cycle, hyperbolic plane, inversion, Malfatti's construction problem

\section{ The history of the problem}

Malfatti (see \cite{malfatti}) raised and solved the following problem:
\emph{construct three circles into a triangle so that each of them touches the two others from outside moreover touches two sides of the triangle too.}
Malfatti determined the radii of the required circles giving an analytic-geometric solution. This formula for the radii was only the beginning of a research motivated by the original problem. We would like to mention here some nice further investigations.

The first nice moment was Steiner's construction. He gave an elegant method (without proof) to construct the given circles. He also extended the problem and his construction to the case of three given circles instead of the sides of a triangle (see in \cite{steiner 1}, \cite{steiner 2}). Cayley referred to this problem in \cite{cayley} as \emph{Steiner's extension of Malfatti's problem}.
We note that Cayley investigated and solved its generalization in \cite{cayley}, he called it also Steiner's extension of Malfatti's problem. His problem is \emph{to determine three conic sections so that each of them touches the two others, and also touches two of three more given conic sections.}
Since the case of circles on the sphere is a generalization of the case of circles of the plane (as it can be seen easily by stereographic projection) Cayley indirectly proved Steiner's second construction. We also have to mention Hart's nice geometric proof for Steiner's construction which was published in \cite{hart}. (It can be found in various textbooks e.g. \cite{casey} and also on the web.)

The second moment which I would like to mention here is two short papers  \cite{schellbach 1} and \cite{schellbach 2} written by Dr.Schellbach. He solved the original construction problem and its "spherical variation", too. The elegant goniometric solution determines the touching points of the circles on the sides of the triangle (spherical triangle). These nice and surprising results inspired Cayley to suggest a little bit more elegant way to get Schellbach's formulas by a change of notation \cite{cayley 2}. Then he supposed that the sides are infinitely small and he reduced the problem to that of a plane triangle. In this way he got some formulas leading to the formulas of his own previous paper \cite{cayley 3}. Finally he expressed the distances from the points of contact of the Malfatti circles to the end-points of the respective sides of the triangle. He noted also that the equations are very similar in form to those given in the same paper for the determination of the radii of the inscribed Malfatti's circles. At this point it seemed that there is not possible to give any new observation for this problem. However Bottema (in \cite{bottema}) observed a simple solution by inversion for the case when the determining figures are touching circles on the euclidean plane. The idea of the construction is that the center of the reference circle of the inversion must be the point of contact of two given circles. Then the inversion leads to a configuration in which the circles can be got in a simple way, and their radii can also be calculated from the radii of the given circles.

It is an interesting fact that nobody has investigated the hyperbolic case of this -- more than two centuries ago raised -- problem before, while the case of the sphere was solved simultaneously with the Euclidean case. Our aim is to compensate for this shortage.

\section{The case of hyperbolic triangle}

Note that the case of the hyperbolic triangle can be solved immediately by the method of Steiner and also by the method of Schellbach. The problem now is the following:  \emph{Construct three circles into a triangle so that each of them touches the two others from outside, and moreover touches two sides of the triangle.}

\subsection{Steiner's construction for hyperbolic triangle}

To solve Malfatti's problem Steiner proposed the following construction \cite{steiner 2} (cf. Figure 2):

\begin{enumerate}
\item Draw three angle bisectors $OA$, $OB$ and $OC$. In the triangles $\triangle_{OAB}$, $\triangle_{OBC}$, $\triangle_{OCA}$ inscribe circles $c'_C$, $c'_A$, $c'_B$, respectively. (The point $O$ is not noted in figure.)
\item For each pair of the circles consider the second (distinct from the corresponding angle bisector) internal tangent lines. These lines concur in a point $K$ and cross the sides in points $H,I$; $D,E$; and $F,G$; respectively. The notations are chosen so that $I\in [B,C]$, $E\in [C,A]$, $G\in [A,B]$; the points $H,D,F$ can lie on any sides different from $[B,C]$, $[C,A]$, $[A,B]$, respectively (where $[X,Y]$ means the segment between the points $X$ and $Y$).
\item The three quadrilaterals $KICE$, $KGBI$, and $KEAG$ have incircles $c_C$, $c_B$, $c_A$, respectively. Their incircles solve Malfatti's problem.
\end{enumerate}

The proof of Hart \cite{hart} works on the hyperbolic plane, too. It is based on two elementary but absolute lemmas (our Lemmas 1 and 2), implying that their corollaries are also valid on the hyperbolic plane. Before these lemmas we recall some concepts.

\begin{figure}[htbp]
\centerline{\includegraphics[scale=1]{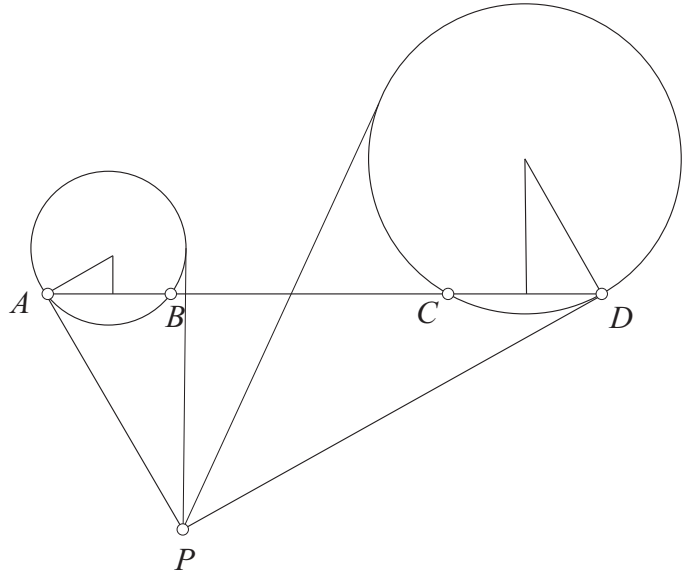}}
\caption{Hart's Lemma 2}
\end{figure}
For two circles we call a \emph{common tangent segment} any segment lying on a common tangent line of the two circles, connecting points of tangency of the two circles. If one of the circles degenerates to a point, then analogously, a \emph{tangent segment between a point and a circle} is any segment lying on a tangent line of the circle, passing through the point, connecting the point and the point of tangency of the circle. Moreover, for points $A,B$, we denote their usual hyperbolic distances by $|AB|$.

\begin{lemma}[\cite{hart}]
A point $P$ is on a common tangent line of two circles if and only if the sum, or the difference of the lengths of two tangent segments, from the point $P$ to the two circles, is equal to the length of some common tangent segment of the two circles.
\end{lemma}

Hence it is evident that if some common tangent lines to each pair of three circles pass through the same point, one of the respective common tangent segments must be equal to the sum or difference of the other two ones, and that the other three common tangent lines (inner or outer ones, as the original ones, respectively) will also pass through a point.

\begin{lemma}[\cite{hart}]
Suppose that two circles cut off segments $AB$ and $CD$ of equal lengths (Fig. 1) from a given line, and the tangent lines to the respective circles at the extreme points $A$, $D$ intersect at $P$. Then the circles will subtend equal angles at $P$. Also if common tangent segments are drawn from each point $A$ and $D$ to the other circle, then they will have equal lengths.
\end{lemma}

The proof of the first lemma is based on two facts; on the triangle inequality and on the fact that the common tangent segments from a point to a circle have equal lengths. Both of them are valid on the hyperbolic plane, too. The proof of the second one uses trigonometry and can be done in the hyperbolic plane, too. Hart himself observed that his proof works on the sphere. Hart's proof of the construction (which we are citing from the original paper because of its elegancy and clarity) is the following (see Fig. 2.):

\begin{figure}[htbp]
\centerline{\includegraphics[scale=1]{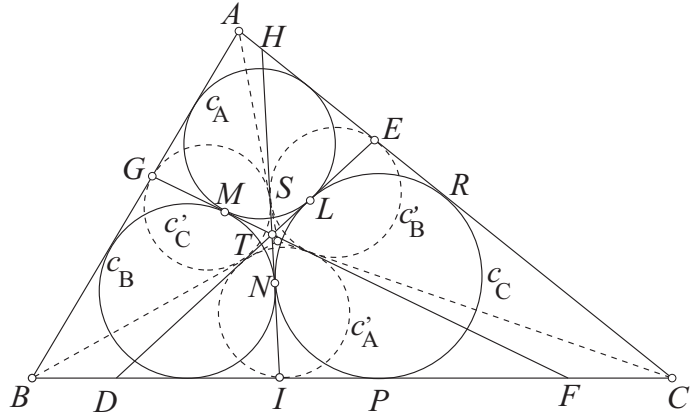}}
\caption{Hart's proof of Steiner's construction.}
\end{figure}

Let $L$, $M$, $N$, be the points of contact of three circles which are pairwise touching, and each touches two sides of the given triangle $\triangle_{ABC}$ (with vertices $A$, $B$, $C$). Draw lines $DE$, $FG$, $HI$, touching these circles at $L,M,N$ and meeting one another at $K$ (observe that they are axes of power). Then since $FI-ID= FO-DP= FM-DL= FK-DK$, $I$ is the point of contact of the circle $c'_{A}$ inscribed in the triangle with the side lines $BC, DE, FG$, with the side $BC$. Similarly, $c'_{B}$ or $c'_{C}$ denote the incircles of the triangles bounded by the lines $CA,FG,HI$, or $AB,HI,DE$, respectively. Let $Q$ be the tangent point of the circle $c'_{A}$ and the line $DE$ (this point is not noted in the figure). Similarly $E$ or $G$ are the points of contact of the circles $c'_{B}$ or $c'_{C}$ with the sides $CA$ or $AB$, respectively. But $IN=IP=QL$, and $NS=ER=EL$, therefore $IS=EQ$, and therefore the circles $c'_{A}$ and $c'_{B}$ subtend equal angles at $C$ (cf. Lemma 2). Also the three common tangents $DE$, $IH$, and $GF$ of the circles $c'_{A}$, $c'_{B}$, $c'_{C}$, pass through $K$; therefore $C$ must be a point on the other common tangent to $c'_{A}$ and $c'_{B}$ (cf. Lemma 1). Similarly, it is proved that the bisectors of the angles $CAB\measuredangle$ or $ABC\measuredangle$ are common tangents to $c'_{B}$ and $c'_{C}$, or to $c'_{C}$ and $c'_{A}$, respectively. Hence Steiner's construction is evident for plane and spherical triangles.

By inversion Hart extends his proof to all configurations of the plane or the sphere in which the sides of the given triangle are arcs of circles. We note that inversion cannot change intersecting circles into non-intersecting ones, thus this method is	
inappropriate to prove Steiner's extension of Malfatti's problem in the case when the given circles are pairwise disjoint.

\subsection{Schellbach's solution}

Schellbach's solution is a goniometric calculation. Historically, this solution preceded Hart's proof (see \cite{schellbach 1} and \cite{schellbach 2} ). To use this for the hyperbolic case we have to change Schellbach's equalities on spherical triangles to the corresponding equalities of hyperbolic triangles. Practically we can use imaginary values to the measure of the lengths and so the trigonometric functions of these quantities give hyperbolic functions depending on their imaginary part. Formally, we follow Cayley's simplified terminology.
Let the sides of the triangle be $a=|BC|$, $b=|CA|$, $c=|AB|$ and let $x,y,z$ be the lengths of the tangent segments from $A,B,C$ to $c_A,c_B,c_C$, respectively (cf. Fig. 2).  Then writing
$$
a+b+c=2s,
$$
$$
a-\frac{1}{2}s=l, \quad b-\frac{1}{2}s=m, \quad c-\frac{1}{2}s=n
$$
whence $l+m+n=\frac{1}{2}s$, and putting also
$$
\frac{1}{2}s-x=\xi, \quad \frac{1}{2}s-y=\eta, \quad \frac{1}{2}s-z=\zeta,
$$
we have
\begin{eqnarray*}
\frac{\cosh l\cosh \eta\cosh \zeta}{\cosh \frac{1}{2}s}+\frac{\sinh l\sinh \eta\sinh \zeta}{\sinh \frac{1}{2}s}& = &1, \\
\frac{\cosh m\cosh \zeta\cosh \xi}{\cosh \frac{1}{2}s}+\frac{\sinh m\sinh \zeta\sinh \xi}{\sinh \frac{1}{2}s}& = &1, \\
\frac{\cosh n\cosh \xi\cosh \eta}{\cosh \frac{1}{2}s}+\frac{\sinh n\sinh \xi\sinh \eta}{\sinh \frac{1}{2}s}& = &1,
\end{eqnarray*}
from which equations the unknown quantities $\xi, \eta, \zeta$ are to be determined. To solve the equations, let the subsidiary angles $\lambda, \mu, \nu$ be determined by the conditions
\begin{eqnarray*}
\frac{\cosh \lambda \cosh m\cosh n}{\cosh \frac{1}{2}s}-\frac{\sinh \lambda \sinh m\sinh n}{\sinh \frac{1}{2}s}& = &1, \\
\frac{\cosh \mu\cosh n\cosh l}{\cosh \frac{1}{2}s}-\frac{\sinh \mu\sinh n\sinh l}{\sinh \frac{1}{2}s}& = &1, \\
\frac{\cosh \nu \cosh l\cosh m}{\cosh \frac{1}{2}s}-\frac{\sinh \nu \sinh l\sinh m}{\sinh \frac{1}{2}s}& = &1.
\end{eqnarray*}
Then it may be shown that
$$
\cosh (\eta +\zeta)=\frac{\cosh \left(\frac{s+\lambda-l}{2}\right)}{\cosh\left(\frac{\lambda+l}{2}\right)}, \quad
\cosh (\eta -\zeta)=\frac{\cosh \left(\frac{s-\lambda+l}{2}\right)}{\cosh\left(\frac{\lambda+l}{2}\right)},
$$
$$
\cosh (\zeta +\xi)=\frac{\cosh \left(\frac{s+\mu-m}{2}\right)}{\cosh\left(\frac{\mu
+m}{2}\right)}, \quad
\cosh (\zeta -\xi)=\frac{\cosh \left(\frac{s-\mu+m}{2}\right)}{\cosh\left(\frac{\mu+m}{2}\right)},
$$
$$
\cosh (\xi +\eta)=\frac{\cosh \left(\frac{s+\nu-n}{2}\right)}{\cosh\left(\frac{\nu+n}{2}\right)}, \quad
\cosh (\xi -\eta)=\frac{\cosh \left(\frac{s-\nu+n}{2}\right)}{\cosh\left(\frac{\nu+n}{2}\right)}.
$$
If we write
\begin{eqnarray*}
\tanh \phi & = & \tanh m\tanh n \coth \frac{1}{2}s, \\
\tanh \chi & = & \tanh n\tanh l \coth \frac{1}{2}s, \\
\tanh \psi & = & \tanh l\tanh m \coth \frac{1}{2}s,
\end{eqnarray*}
then
\begin{eqnarray*}
\cosh(\lambda- \phi)& = & \frac{\cosh \frac{1}{2}s\cosh \phi}{\cosh m}{\cosh n}, \\
\cosh(\mu- \chi)& = & \frac{\cosh \frac{1}{2}s\cosh \chi}{\cosh n}{\cosh l}, \\
\cosh(\nu- \psi)& = & \frac{\cosh \frac{1}{2}s\cosh \psi}{\cosh l}{\cosh m},
\end{eqnarray*}
which give the values of $\lambda, \mu, \nu$. From these quantities $\xi, \eta, \zeta$ and thus $x, y, z$ can be constructed, too.

\section{The case of the cycles}

Our aim is to investigate touching cycles of the hyperbolic plane. For simplicity, we assume that the interiors of the given cycles are pairwise non-overlapping. The interior of a paracycle in the Poincar\'e model is the interior of the circle representing it, and the interior of a hypercycle is the locus of the points with distances from its base line smaller than
the constant used in the definition of the hypercycle.
We formulate the exercise as follows: \emph{Determine three cycles of the hyperbolic plane so that each of them touches the two others moreover touches two of three given cycles of the hyperbolic plane.}

\begin{figure}[htbp]
\centerline{\includegraphics[scale=1]{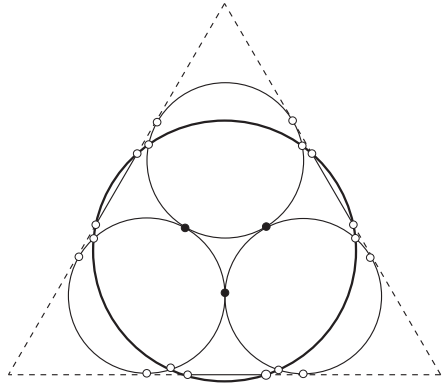}}
\caption{The euclidean solution does not work on the hyperbolic plane.}
\end{figure}

Using Poincar\'e's conformal disk model the given cycles are represented by the arcs of certain circles. By Steiner's method we can construct the corresponding Malfatti's circles solving Steiner's extension of Malfatti's problem. But it is possible that the corresponding arcs of circles have no common point in the model-circle. This means that there is no solution of the problem with respect to the hyperbolic plane. In fact, consider the situation of Fig.3. The edges of the regular triangle with respect to the basic circle of the model represent arcs of hypercycles of the hyperbolic plane. On the other hand the uniquely determined Malfatti's circles with respect to this regular triangle touch the corresponding edges in external points of the model implying that the corresponding cycles cannot be touching cycles to the given hypercycles.

The trouble of this example can be solved easily if we consider the whole cycle (cf. Def.1) rather than one of its connected components. We will use the following concepts:

\begin{defi}
The \emph{hypercycle} is the locus of points of the hyperbolic plane with distances from a line (called its \emph{base line}) being equal to a constant. (Its interior is convex, and is bounded by a curve with two connected components.) \emph{Two cycles are touching} if they have a common point with a common tangent line.
\end{defi}

We can state the following theorem on existence.

\begin{theorem}
For three given cycles there are three other cycles such that each of them touches the two others moreover touches two from the three given ones.
\end{theorem}

Remark that in this theorem the interiors of the corresponding six cycles need not to be disjoint.  On the other hand the tools of a concrete construction will require some further geometric conditions thus in the rest of this paper we assume that the interiors of the given cycles are pairwise non-overlapping.

\vskip0.5cm
\proof Consider the Poincar\'e's disk model of the hyperbolic plane. In the case, when there is no hypercycle among the given cycles the statement is trivial. On the other hand by our definition a hypercycle is the union of two arcs of circles in this model. Now consider for a hypercycle the entire circle containing one of the original arcs. The inner arc (which lies in the interior of the model) of the chosen circle is the inverse image (with respect to the model circle) of the outer arc of the other possible circle. Now Steiner's construction can be done on the embedding Euclidean plane and we get a Malfatti's system of circles for the representing circles. If we take the inverse of those points of the six circles which are external to the model, we get either whole hypercycle; paracycle; or a circle of the hyperbolic plane accordingly to the cases when the respective circle intersects; touches externally; or lies in the complement of the model circle, respectively.

So with this transformation from a Malfatti's system of circles of the embedding Euclidean plane, we get a Malfatti's system of cycles of the hyperbolic plane corresponding to the given cycles as we stated.
\qed

Unfortunately the used condition for touching cycles cannot exclude the possibility of existence of other common points of the two cycles (contrary to the case of circles).
Two possible points of intersection can lie on the other connected components of the hypercycle. (On the example of Fig.3 we can find this situation after the inversion of the external parts of the hypercycles.)

Thus with more rigorous definitions of touching we have new (typically harder) problems on existence.  Two possibilities for the concepts of touching are:

\begin{itemize}
\item two cycles are touching if they have exactly one common (real) point on the hyperbolic plane (in that point their tangent lines are also the same), or
\item two cycles are touching if their interiors are touching externally.
\end{itemize}

In this paper we do not deal with these questions which can give the base of a foregoing paper. We consider here the most simple case of touching as we fixed in our Definition 1. Our question is:

\begin{center}
\emph{How can we construct a Malfatti system of cycles of a given system of cycles?}
\end{center}

\subsection{Steiner's construction on the hyperbolic plane}

In this section we present that possible form of Steiner's construction which best meet the exception of the original problem. We note (see the discussion in the proof) that Theorem 2 has a more general form giving all possible solutions of the problem, however for simplicity we restrict ourself to the most plausible case, when the cycles touch each other from outside. Throughout this section we use the fact that cycles represented by circles in the conformal model of Poincar\'e. The Euclidean constructions on circles of this model gives hyperbolic constructions on cycles in the hyperbolic plane. To do these constructions manually we have to use special rulers and calipers to draw the distinct types of cycles. For brevity, we think for a fixed conformal model of the embedding Euclidean plane and preserve the name of the known Euclidean concepts with respect to the corresponding concept of the hyperbolic plane, too.

\begin{theorem}
Steiner's construction can be done also in the hyperbolic plane. More precisely, for three given non-overlapping cycles there can be constructed three other cycles, each of them touches the two other ones from outside and also touches two of the three given cycles from outside.
\end{theorem}

\proof Denote by $c_i$ the given cycles. Now the steps of Steiner's construction are the following. (These steps are some analogues of steps $1,2,3$ in 1.2.1.)
\begin{enumerate}
\item Construct the cycle of inversion $c_{i,j}$, for the given cycles $c_i$ and $c_j$, where the center of inversion is the external centre of similitude of them. (I.e., the center of $c_{i,j}$ is the center of the above inversion, and $c_i,c_j$ are images of each other with respect to inversion with respect to $c_{ij}$. Observe that $c_{ij}$ separates $c_i$ and $c_j$.)
\item Construct cycle $k_j$ touching two cycles $c_{i,j}, c_{j,k}$ and the given cycle $c_j$, in such a way that $k_j,c_j$ touch from outside, and $k_{ij},c_{ij}$ (or $c_{jk}$) touch in such a way that $k_j$ lies on that side of $c_{ij}$ (or $c_{ik}$) on which side of them $c_j$ lies.
\item Construct the cycle $l_{i,j}$ touching $k_i$ and $k_j$ through the point $P_k=k_k\cap c_k$.
\item Construct Malfatti's cycle $m_j$ as the common touching cycle of the four cycles $l_{i,j}$, $l_{j,k}$, $c_i$, $c_k$.
\end{enumerate}
The first step is the hyperbolic interpretation of the analogous well-known Euclidean construction on circles.

To the second step we follow Gergonne's construction (see e.g. \cite{dorrie}), see in Fig. 4. This construction is:

\emph{Draw the point $P$ of power of the given circles $c_1,c_2,c_3$ and an axis of similitude of certain three centres of similitude.  We join the poles $P_1,P_2,P_3$ of this axis of similitude with respect to the circles $c_1,c_2,c_3$ with the point $P$ by straight lines. Then the lines $PP_i$ cut the circles $c_i$ in two points $Q_{i1}$ and $Q_{i2}$. A suitable choice $Q_{1j(1)},Q_{2j(2)},Q_{3j(3)}$ will give the touching points of some sought circle and $c_1,c_2,c_3$. More exactly, there are two such choices $Q_{1j(1)},Q_{2j(2)},Q_{3j(3)}$ and $Q_{1k(1)},Q_{2k(2)},Q_{3k(3)}$, satisfying $j(i) \ne k(i)$ for $1 \le i \le 3$, where $|PP_{ij(i)}| \le |PP_{ik(i)}|$.}

\begin{figure}[htbp]
\centerline{\includegraphics[scale=1]{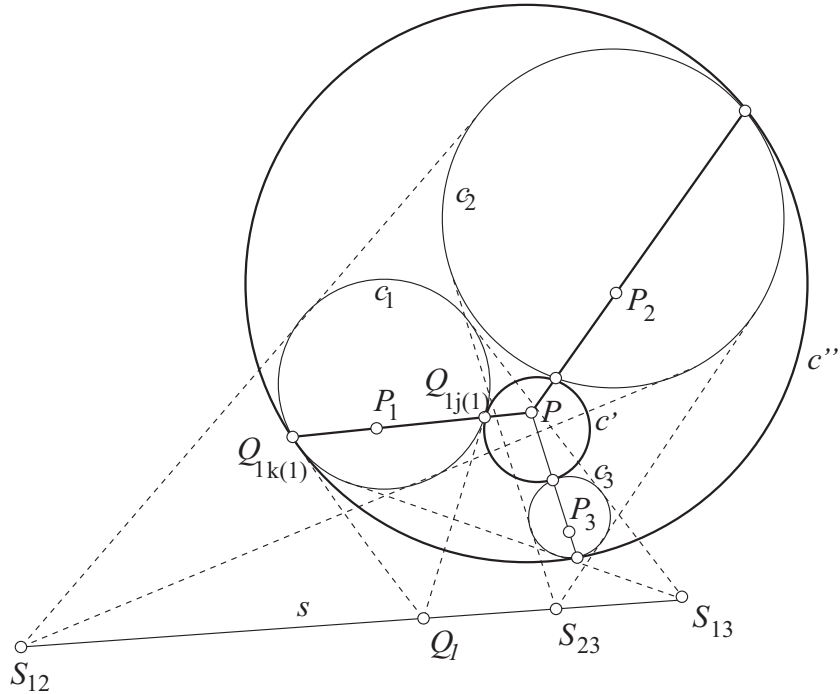}}
\caption{The construction of Gergonne}
\end{figure}

In Fig. 4 the axis of similitude contains the three outer centers of similitude, in which case, choosing for $Q_{ij(i)}$ the intersection points closer to $P$, we obtain the common outward touching cycle, and for choosing the farther intersection points to $P$ we obtain the common touching cycle that contains $c_1,c_2,c_3$. We denoted these circles in Fig. 4 by $c'$ and $c''$, respectively.

Choosing, e.g., for $c_1,c_3$ and $c_2,c_3$ the inner centers of similitude, and then for $c_1,c_2$ the outer center of similitude, we obtain another axis of similitude (by permuting the indices we obtain still two more similar cases). Then defining the points $P_i$ and $P_{ij(i)}$ analogously like above, if $|PQ_{1j(1)}| \le |PQ_{1k(1)}|$, $|PQ_{2j(2)}| \leq |PQ_{2k(2)}|$, and $|PQ_{3j(3)}| \geq |PQ_{3k(3)}|$,
then the circle $Q_{1j(1)}Q_{2j(2)}Q_{3j(3)}$ touches $c_1,c_2,c_3$, contains $c_3$ and touches $c_1,c_2$ externally, while the circle \newline \hfill
$Q_{1k(1)}Q_{2k(2)}Q_{3k(3)}$ touches $c_1,c_2,c_3$, contains $c_1,c_2$, and touches $c_3$ externally.

Summing up: there are eight cycles touching each of $c_1,c_2,c_3$.

An Euclidean proof of the pertinence of this construction on circles can be rewritten also by hyperbolic terminology. Consider the cycles $c'$ and $c''$ touching $c_1$, $c_2$ and $c_3$, in any of the four above described cases; in Fig. 4 $c'$ touches each of $c_1,c_2,c_3$ externally, and $c''$ touches each of $c_1,c_2,c_3$ internally. Then the line joining the touching points $Q_{ij(i)}$ and $Q_{ik(i)}$ passes through one of the centers of similitude $P$ of $c'$ and $c''$. Thus $P$ is the point of power of $c_1$, $c_2$ and $c_3$. On the other hand, two of the three given cycles (say $c_1$ and $c_2$) give a touching pair with respect to $c'$ and $c''$, hence its outer center of similitude $S_{12}$ has the same power with respect to $c'$ and $c''$. So the three outer centers of similitude $S_{12}$, $S_{13}$ and $S_{23}$ are on the axis of power of $c'$ and $c''$. (It is also (by definition) an axis of similitude with respect to $c_1$, $c_2$ and $c_3$, say $s$. For $c',c''$ being another pair of touching circles, in the other three cases, the respective changes have to be made in the choice.) Since the pole $Q_i$ (with respect to the cycle $c_i$) of the line joining $Q_{ij(i)}$ and $Q_{ik(i)}$ is the intersection point of the common tangents of $c'$ and $c_i$ at $Q_{ij(i)}$, and  $c''$ and $c_i$ at $Q_{ik(i)}$, respectively, it is also on $s$. By the theorem of pole-polar we get that the pole $P_i$ of $s$ with respect to $c_i$ lies on the line $Q_{ij(i)}Q_{ik(i)}$. This proves the construction.

The third step is a special case of the second one. (A given cycle is a point now.) Obviously the general construction can be done in this case, too.

The fourth step is again the second one choosing three arbitrary cycles from the four ones since the quadrangles determined by the cycles have incircles.

Finally we have to prove that this construction gives the Malfatti's cycles. As we saw the Malfatti's cycles are exist (see Theorem 1). We also know that in an embedding hyperbolic space the examined plane can be inverted to a sphere. The trigonometry of the sphere is absolute implying that the possibility of a construction which can be checked by trigonometric calculations, is independent of the fact that the embedding space is a hyperbolic space or a Euclidean one. Of course, the Steiner construction is just such a construction, the touching position of circles on the sphere can be checked by spherical trigonometry. So we may assume that the examined sphere is a sphere of the Euclidean space and we can apply Cayley's analytical research (see in \cite{cayley}) in which he proved that Steiner's construction works on a surface of second order. Hence the above construction produces the required touchings.
\qed

\begin{remark}
We note that all the notions and constructions in the proof Theorem 2 can be discussed without the use of a model. A part of the basis of this ``absolute'' description is given in C. V\"or\"os
(cf. \cite{voros}), who  introduced \emph{general points} (which can be real or ideal points, or points at infinity), and defined the hyperbolic length of a segment with general endpoints.
Combining it with Casey's approach on spherical geometry (cf. \cite{casey 1}), a model independent proof can be obtained. Nevertheless, in light of the length and complexity of this process, we will publish it  in a forthcoming paper.
\end{remark}

\section{Acknowledgement}
I thank for the rigorous criticism and valuable suggestions of the referee, he helped me to write a more readable and understandable form this paper.

\begin{center}
\'Akos G.Horv\'ath \\  Department of Geometry, Mathematical Institute \\
Budapest University of Technology and Economics,\\
H-1521 Budapest,\\
Hungary\\
e-mail: ghorvath@math.bme.hu
\end{center}

\end{document}